\documentclass[a4paper,11pt]{amsart}

\usepackage{graphicx}
\usepackage{mathptmx}
\usepackage{amsmath}
\usepackage{amssymb}
\usepackage{enumitem}
\usepackage{xcolor}
\usepackage{pgfplots}

\newmuskip\pFqmuskip

\newcommand*\pFq[6][8]{%
  \begingroup 
  \pFqmuskip=#1mu\relax
  \mathcode`=\string"8000
  \begingroup\lccode`\~=`\,
  \lowercase{\endgroup\let~}\pFqcomma
  F^{#2}_{#3}{\left(\genfrac..{0pt}{}{#4}{#5}\bigg|#6\right)}%
  \endgroup
}
\newcommand{\pFqcomma}{\mskip\pFqmuskip}

\newtheorem{theorem}{Theorem}[section]

\begin{document}

\title[Recurrence relations for harmonic and derangement numbers]{Recurrence relations for harmonic and derangement numbers}

\author{Taekyun  Kim}
\address{Department of Mathematics, Kwangwoon University, Seoul 01897, Republic of Korea}
\email{tkkim@kw.ac.kr}
\author{Dae San  Kim}
\address{Department of Mathematics, Sogang University, Seoul 04107, Republic of Korea}
\email{dskim@sogang.ac.kr}

\author{JongKyum Kwon}
\address{Department of Mathematics Education and ERI, Gyeongsang National University, Jinju 52828, Republic of Korea}
\email{mathkjk26@gnu.ac.kr }

\author{Kyo-Shin Hwang}
\address{Graduate School of Education, Yeungnam University, Gyeongsan 38541, Republic of Korea}
\email{kshwang@yu.ac.kr}


\begin{abstract}
We use elementary methods to establish three key recurrence relations: one for derangement numbers, a second for harmonic numbers, and a third for degenerate harmonic numbers. Our results not only contribute to the understanding of the underlying structure of these numbers but also highlight the effectiveness of elementary techniques in discovering new mathematical properties. The findings have potential applications in various fields where these numbers appear, including combinatorics, probability, and computer science. 
\end{abstract}

\maketitle

\section{Introduction}
In this paper, we show by using elementary method the following recurrence relations for the derangement numbers $D_{n}$, the harmonic numbers $H_{n}$, and the degenerate harmonic numbers $H_{n,\lambda}$: for any integers $m,n \ge 0$, we have
\begin{align*}
&\frac1{n!}{D_{m+n}}=\sum_{l=0}^n  \sum_{k=0}^m \binom{k+n-l-1}{n-l}\binom{m}{k} (-1)^{m-k}\frac{k!}{l!} D_l, \label{0} \\
& \binom{m+n}{n}H_{n+m}= \sum_{k=0}^n H_k  \binom{m+n-k-1}{n-k} + H_m \binom{m+n}{n}, \nonumber \\
&\binom{m+n}{n} H_{n+m,\lambda}= \sum_{k=0}^n  H_{k,\lambda} \binom{m+n-k-1}{n-k}  + H_{m,\lambda} \binom{m+n-\lambda}{n}. \nonumber
\end{align*} \par
The exploration of degenerate versions of special polynomials and numbers (see [1,6-10,16] and the references therein) began with Carlitz's foundational work on degenerate Bernoulli and Euler numbers. This field has since expanded to include not only polynomials and numbers but also transcendental functions like the gamma function. In addition, the development of the $\lambda$-umbral calculus has provided a more robust and effective framework for analyzing degenerate Sheffer polynomials compared to the classical umbral calculus. \par
A derangement is a permutation of a set of elements where no element remains in its original position. The number of derangements for a set of $n$ elements is known as the $n$-th derangement number, denoted as $D_{n}$. A classic and intuitive example of a derangement is the hat-check problem: Imagine $n$ guests at a party each check their hat, and the hats are returned randomly. The number of ways that no guest receives their own hat is the derangement number, $D_{n}$. \par
Harmonic numbers $H_{n}$ have a rich history and are found throughout various fields of mathematics, computer science, probability, physics, and engineering. Their widespread appearance underscores their fundamental importance in many different areas of study. The degenerate harmonic numbers $H_{n,\lambda}$ are a degenerate version of the harmonic numbers. \par

We have structured this paper into three sections. Section 1 serves as an introduction to the core concepts used throughout the paper. We begin by recalling derangement and harmonic numbers, followed by a review of degenerate exponentials and logarithms. This section also introduces the notions of degenerate harmonic numbers and hyperharmonic numbers. The main contributions of this work are found in Section 2, where we establish several important recurrence relations. Specifically, Theorem 2.1 presents a recurrence for derangement numbers, Theorem 2.2 for harmonic numbers, and Theorem 2.3 for degenerate harmonic numbers. Finally, Section 3 provides a summary of our results and concludes the paper. A list of general references can be found in [3,4,12,13]. In the rest of this section, we recall the facts that are needed throughout this paper.\par
\vspace{0.1in}
The $n$-the derangement $D_{n}$ is given by
\begin{equation}\label{1}
D_n = n! \sum_{k=0}^n \frac{(-1)^k}{k!}, \quad(\rm see, \ [3,9]),
\end{equation}
where $n$ is a nonnegative integer. From \eqref{1}, we have
\begin{equation}\label{2}
\frac1{1-t}e^{-t} =\sum_{n=0}^\infty D_n \frac{t^n}{n!},\quad(\rm see, \ [3,9]).
\end{equation}
The harmonic numbers are defined by (see [2,5,11,14,16,17])
\begin{equation}\label{3}
H_0=0,\quad H_n = 1 + \frac12 + \cdots +\frac1n, \quad(n \ge 1).
\end{equation}
From \eqref{3}, we have (see [7,8,15])
\begin{equation}\label{4}
\frac1{1-t}\log\left(\frac1{1-t}\right)=\sum_{n=1}^\infty H_n t^n.
\end{equation} \par
For any nonzero $\lambda\in \mathbb{R}$, the degenerate exponentials are defined by (see [6-9,15])
\begin{equation*}
e_\lambda^x (t)= \sum_{n=0}^\infty (x)_{n,\lambda}\frac{t^n}{n!},\quad e_\lambda =e_\lambda^{1}(t),
\end{equation*}
where
$$(x)_{0,\lambda}=1, \quad (x)_{n,\lambda}=x(x-\lambda)\cdots(x-(n-1)\lambda), \quad(n\ge 1).$$
The degenerate logarithm $\log_{\lambda}(t)$ is defined as the compositional inverse of $e_\lambda(t)$ and given by
\begin{equation}\label{5}
\log_\lambda(t)= \frac{1}{\lambda}(t^{\lambda}-1).
\end{equation}
Note here that $\lim_{\lambda \rightarrow 0}\log_{\lambda}(t)=\log(t)$.
We use the following relation of the degenerate logarithms:
\begin{equation}\label{6}
\log_{\lambda}(xy)=\log_{\lambda}(x)+x^{\lambda}\log_{\lambda}(y)=\log_{\lambda}(y)+y^{\lambda}\log_{\lambda}(x).
\end{equation}
We note from \eqref{5} that
\begin{equation*}
\log_\lambda(1+t)=\sum_{n=1}^\infty \lambda^{n-1}(1)_{n,1/\lambda}\frac{t^n}{n!}
=\sum_{n=1}^\infty \binom{\lambda-1}{n-1}\frac{t^n}{n},\quad(\rm see\ [7,8,9]).
\end{equation*} \par
Kim-Kim defined the degenerate harmonic numbers given by (see [7,8,15])
\begin{equation}\label{7}
H_{0,\lambda}=0,\quad H_{n,\lambda}= \frac1{\lambda}\sum_{k=1}^n\binom{\lambda}{k}(-1)^{k-1} =
\sum_{k=1}^n\binom{\lambda-1}{k-1}\frac{(-1)^{k-1}}{k}, \quad(n\ge 1).
\end{equation}
Note that $$\lim_{\lambda\to 0} H_{n,\lambda}=H_n,\quad(n\ge 0).$$
From \eqref{7}, we derive the generating function of degenerate harmonic numbers:
\begin{equation}\label{8}
\frac1{1-t}\log_{-\lambda}\left(\frac1{1-t}\right)=\sum_{n=1}^\infty H_{n,\lambda}t^n, \quad H_{0,\lambda}=0,\quad(\rm see\ [7,8,15]).
\end{equation}

In 1996, Conway and Guy introduced the hyperharmonic numbers, $H_{n}^{(r)},\ (n,r\ge 0),$ which are defined by (see \cite{4})
\begin{equation}\label{9}
H_{0}^{(r)}=0,\,\,(r \ge 0),\quad H_{n}^{(0)}=\frac{1}{n},\,\,(n \ge 1),\quad H_{n}^{(r)}=\sum_{k=1}^{n}H_{k}^{(r-1)},\,\, (n,r\ge 1).
\end{equation}
Note that $H_{n}^{(1)}=H_{n}$ are the harmonic numbers.
From \eqref{9}, we have (see \cite{4})
\begin{equation}\label{10}
H_{n}^{(m+1)}=\binom{n+m}{m}\big(H_{n+m}-H_{m}\big),\ (m\ge 0), 	
\end{equation}
and
\begin{equation}\label{11}
\frac{1}{(1-t)^{r}} \log\left(\frac{1}{1-t}\right)=\sum_{n=1}^{\infty}H_{n}^{(r)}t^{n}.
\end{equation} \par
Recently, Kim-Kim defined the degenerate hyperharmonic numbers, $H_{n,\lambda}^{(r)},\ (n,r\ge 0)$, which are given by (see [7,8,15])
\begin{align}\label{12}
&H_{0,\lambda}^{(r)}=0,\ H_{n,\lambda}^{(0)}=\frac{1}{n!}\lambda^{n-1}(-1)^{n-1}(1)_{n,\frac{1}{\lambda}}=\frac{1}{\lambda}\binom{\lambda}{n}(-1)^{n-1}, \ (n \ge 1), \\
&H_{n,\lambda}^{(r)}=\sum_{k=1}^{n}H_{k,\lambda}^{(r-1)},\ (n,r \ge 1). \nonumber
\end{align}
Observe that $H_{n,\lambda}^{(1)}=H_{n,\lambda}$ are the degenerate harmonic numbers. From \eqref{12}, we note that (see [7])
\begin{equation}\label{13}
(-1)^{m}\binom{\lambda-1}{m}H_{n,\lambda}^{(m+1)}=\binom{n+m}{m}(H_{n+m,\lambda}-H_{m,\lambda}),\ (m\ge 0), 
\end{equation}
and (see [8])
\begin{equation*}
\frac{1}{(1-t)^{r}}\log_{-\lambda}\left(\frac{1}{1-t}\right)=\sum_{n=1}^{\infty}H_{n,\lambda}^{(r)}t^{n},\quad (r\ge 0). 	
\end{equation*}
Note that (see \eqref{11})
\begin{displaymath}
\lim_{\lambda\rightarrow 0}H_{n,\lambda}^{(r)}=H_{n}^{(r)},\quad (n\ge 0).
\end{displaymath} \par

\section{Recurrence relations for harmonic and derangement numbers}
In the same spirit as [6], we derive three recurrence relations for the derangement numbers $D_{n}$, the harmonic numbers $H_{n}$, and the degenerate harmonic numbers $H_{n,\lambda}$.
From \eqref{2}, we note that
\begin{equation}\label{14}
\begin{split}
&\sum_{m,n=0}^\infty D_{n+m}\frac{x^n}{n!}\frac{y^m}{m!}=
\sum_{m=0}^\infty\frac{d^m}{dx^m}\sum_{n=0}^\infty D_n \frac{x^n}{n!}\frac{y^m}{m!}\\
&=\sum_{n=0}^\infty\frac{D_n}{n!}\sum_{m=0}^n \binom{n}{m} {x^{n-m}}{y^m}
=\sum_{n=0}^\infty\frac{D_n}{n!}(x+y)^n\\
&=\frac1{1-x-y}e^{-(x+y)} = \frac{e^{-x}}{1-x} \frac1{1-\frac{y}{1-x}} e^{-y}\\
&=\frac{e^{-x}}{1-x}\sum_{k=0}^\infty\left(\frac1{1-x}\right)^k y^k \sum_{m=0}^\infty\frac{(-1)^m}{m!}y^m\\
&=\frac1{1-x}e^{-x}\sum_{k=0}^\infty\left(\frac1{1-x}\right)^k y^k k! \sum_{m=k}^\infty
\frac{(-1)^{m-k}m!}{(m-k)!k!} \frac{y^{m-k}}{m!}\\
&=\frac1{1-x}e^{-x}\sum_{k=0}^\infty k! \left(\frac1{1-x}\right)^k \sum_{m=k}^\infty
(-1)^{m-k}\binom{m}{k} \frac{y^{m}}{m!}\\
&=\sum_{m=0}^\infty \frac{y^m}{m!} \sum_{k=0}^m k! \left(\frac1{1-x}\right)^k
(-1)^{m-k}\binom{m}{k}  \frac1{1-x}e^{-x}\\
&=\sum_{m=0}^\infty \frac{y^m}{m!} \sum_{k=0}^m k! \sum_{j=0}^\infty \binom{k+j-1}{j} x^j
(-1)^{m-k}\binom{m}{k} \sum_{l=0}^\infty D_l  \frac{x^l}{l!}\\
&=\sum_{m=0}^\infty  \sum_{n=0}^\infty \frac{y^m}{m!} \frac{x^n}{n!}\left(n!\sum_{l=0}^n
 \sum_{k=0}^m \binom{k+n-l-1}{n-l}\binom{m}{k}
(-1)^{m-k}\frac{k!}{l!} D_l \right).
\end{split}
\end{equation}
Therefore, by comparing the coefficients on both sides of \eqref{14}, we obtain the following theorem.
\par
\begin{theorem}
For any integers $m, n\ge 0$, we have
$$\frac1{n!}{D_{m+n}}=\sum_{l=0}^n  \sum_{k=0}^m \binom{k+n-l-1}{n-l}\binom{m}{k} (-1)^{m-k}\frac{k!}{l!} D_l.$$
\end{theorem}
From \eqref{4}, we have
\begin{equation}\label{15}
\begin{split}
&\sum_{m,n=0}^\infty\binom{m+n}{m}H_{n+m}x^ny^m
=\sum_{m=0}^\infty\frac{d^m}{dx^m}\sum_{n=0}^\infty H_n x^n \frac{y^m}{m!}\\
&=\sum_{n=0}^\infty H_n \frac{d^m}{dx^m} \sum_{m=0}^\infty  x^n \frac{y^m}{m!}
=\sum_{n=0}^\infty H_n \sum_{m=0}^n  \binom{n}{m}   x^{n-m}{y^m}\\
&=\sum_{n=0}^\infty H_n \left(x+y\right)^n =\frac1{1-x-y}
\log\left(\frac1{1-x-y}\right)\\
&=\frac1{1-x}\frac1{1-\frac{y}{1-x}}\log\left(\frac1{1-x} \frac1{1-\frac{y}{1-x}}\right)
\end{split}
\end{equation}
\begin{align*}
&=\left(\frac1{1-x}\log\left(\frac1{1-x}\right)\right) \left(\frac1{1-\frac{y}{1-x}}\right) +
\frac1{1-x} \frac1{1-\frac{y}{1-x}} \log\left(\frac1{1-\frac{y}{1-x}}\right)\\
&= \sum_{k=0}^\infty H_k x^k \sum_{m=0}^\infty \left( \frac{y}{1-x}\right)^m + \sum_{m=0}^\infty H_m y^m \left( \frac{1}{1-x}\right)^{m+1}\\
&= \sum_{k=0}^\infty H_k x^k \sum_{m=0}^\infty y^m \sum_{j=0}^\infty \binom{m+j-1}{j} x^j
+ \sum_{m=0}^\infty H_m y^m \sum_{n=0}^\infty\binom{n+m}{n}x^n\\
&=\sum_{m=0}^\infty \sum_{n=0}^\infty     x^n y^m \left( \sum_{k=0}^n H_k  \binom{m+n-k-1}{n-k}
+ H_m \binom{n+m}{n}\right).
\end{align*}
Therefore, by comparing the coefficients on both sides of \eqref{15}, we obtain the following theorem.\par

\begin{theorem}
For any integers $m,n\ge 0$, we have
$$ \binom{m+n}{m}H_{n+m}= \sum_{k=0}^n H_k  \binom{m+n-k-1}{n-k} + H_m \binom{n+m}{n}.$$
\end{theorem}

By \eqref{10} and Theorem 2.2, we get
$$ H_n^{(m+1)}=\binom{n+m}{m}\left(H_{n+m}-H_m \right)= \sum_{k=0}^n H_k \binom{m+n-k-1}{n-k}.$$
Now, we obtain that
\begin{equation}\label{16}
\begin{split}
&\sum_{m,n=0}^\infty H_{n+m,\lambda}\binom{n+m}{m}x^ny^m
=\sum_{m=0}^\infty\frac{d^m}{dx^m}\sum_{n=0}^\infty H_{n,\lambda}x^n\frac{y^m}{m!}\\
&=\sum_{n=0}^\infty H_{n,\lambda}  \sum_{m=0}^\infty    \frac{d^m}{dx^m} x^n\frac{y^m}{m!}
= \sum_{n=0}^\infty H_{n,\lambda}  \sum_{m=0}^n \binom{n}{m}  x^{n-m}{y^m}\\
&=\sum_{n=0}^\infty H_{n,\lambda}(x+y)^n = \frac1{1-x-y}\log_{-\lambda}\left(\frac1{1-x-y}\right).
\end{split}
\end{equation}
On the other hand, by \eqref{6} and \eqref{8}, we get
\begin{equation}\label{17}
\begin{split}
&  \frac1{1-x-y}\log_{-\lambda}\left(\frac1{1-x-y}\right)
= \frac1{1-x} \frac1{1-\frac{y}{1-x}}\log_{-\lambda}\left(\frac1{1-x} \frac1{1-\frac{y}{1-x}}\right)\\
&=\frac1{1-x} \frac1{1-\frac{y}{1-x}}\left(\log_{-\lambda}\left(\frac1{1-x}\right) +
\left(\frac1{1-x}\right)^{-\lambda} \log_{-\lambda}\left(\frac1{1-\frac{y}{1-x}}\right)\right) \\
&=\frac1{1-x} \log_{-\lambda}\left(\frac1{1-x}\right) \frac1{1-\frac{y}{1-x}}
 + \left(\frac1{1-x}\right)^{1-\lambda} \frac1{1-\frac{y}{1-x}}  \log_{-\lambda}\left(\frac1{1-\frac{y}{1-x}}\right)\\
&=\sum_{l=0}^\infty H_{l,\lambda} x^l \sum_{m=0}^\infty y^m \left(\frac1{1-x}\right)^m
+ \left(\frac1{1-x}\right)^{1-\lambda}  \sum_{m=0}^\infty H_{m,\lambda}\left(\frac{y}{1-x}\right)^m\\
&=\sum_{l=0}^\infty H_{l,\lambda} x^l \sum_{m=0}^\infty y^m
\sum_{k=0}^\infty \binom{m+k-1}{k}x^k + \sum_{m=0}^\infty H_{m,\lambda}y^m \left(\frac1{1-x}\right)^{m+1-\lambda}\\
&=\sum_{m=0}^\infty  \sum_{n=0}^\infty y^m x^n \left(\sum_{l=0}^n  H_{l,\lambda}
\binom{m+n-l-1}{n-l}  + H_{m,\lambda} \binom{m+n-\lambda}{n}\right).
\end{split}
\end{equation}
Therefore, by \eqref{16} and \eqref{17}, we obtain the following theorem.
\par

\begin{theorem}
For any integers $m,n\ge 0$, we have
\begin{equation*}
\binom{n+m}{n} H_{n+m,\lambda}= \sum_{l=0}^n  H_{l,\lambda} \binom{m+n-l-1}{n-l}  + H_{m,\lambda} \binom{m+n-\lambda}{n}.
\end{equation*}
\end{theorem}
Thus, by \eqref{13} and Theorem 2.3, we get
\begin{equation*}
H_{n,\lambda}^{(m+1)}= \frac{(-1)^m}{\binom{\lambda-1}{m}}\left[
\sum_{l=0}^n  H_{l,\lambda}\binom{m+n-l-1}{n-l}  + \binom{m+n-\lambda}{n} H_{m,\lambda} -\binom{n+m}{n} H_{m,\lambda}\right].
\end{equation*}
\par
\section{Conclusion}
In this paper, we have successfully derived recurrence relations for the derangement numbers $D_{n}$, the harmonic numbers $H_{n}$, and the degenerate harmonic numbers $H_{n,\lambda}$ using elementary methods. The recurrence relations, as proven in Theorem 2.1, Theorem 2.2, and Theorem 2.3, provide a new and structured way to compute and analyze these numbers. The simplicity of our approach suggests that similar methods could be used to discover new recurrences of other special numbers and polynomials. This research not only enriches the theory of number sequences but also serves as a testament to the power of elementary techniques in uncovering profound mathematical truths. \par
It is one of our future projects to continue to explore special numbers and polynomials by employing various methods.

\end{document}